\documentclass{amsproc}

\usepackage{graphicx}

\theoremstyle{definition}

\theoremstyle{remark}

\numberwithin{equation}{section}

\newcommand{\bt}{{\bf 10}}
\newcommand{\bfv}{{\bf 5}}
\newcommand{\bfb}{{\overline{\bf 5 \!}\,}}
\newcommand{\btb}{{\overline{\bf 10 \!}\,}}
\newcommand{\newsubsection}[1]{\begin{flushleft} \subsection{#1}
\end{flushright}}

\begin{document}

\title{Higgs Bundles and String Phenomenology}

\author{Martijn Wijnholt}

\address{Ludwig-Maximilians-Universit\"at, Arnold-Sommerfeld Center,
Theresienstrasse 37, D-80333 Munich, Germany}


\begin{abstract}
String phenomenology is the branch of string theory concerned with
making contact with particle physics. The original models involved
compactifying ten-dimensional supersymmetric Yang-Mills theory on an
internal Calabi-Yau three-fold. In recent years, this picture has
been extended to compactifications of supersymmetric Yang-Mills
theory in seven, eight or nine dimensions. We review some of this
progress and explain the fundamental role played by Higgs bundles in
this story.

\end{abstract}

\maketitle

\vspace{-.5cm}

\setcounter{tocdepth}{2}
\tableofcontents

\baselineskip=12pt
\parskip 5pt plus 1pt

\vspace{-.7cm}

\section{String Phenomenology}

\subsection{Particle physics and Grand Unification}

$ $

\noindent To date the most successful description of particle
physics is the Standard Model. It is a gauge theory with gauge group
$SU(3)_c \times SU(2)_w \times U(1)_Y$ and with three fermionic
generations of matter:
\begin{equation}
\begin{array}{lll}
  Q = ({\bf 3}, {\bf 2})_{1/6} & u^c = ({\overline{\bf 3
\!}\,},{\bf 1})_{-2/3}& d^c = ({\overline{\bf 3 \!}\,},{\bf
1})_{1/3} \\[2mm]
  L = ({\bf 1}, {\bf 2})_{-1/2} & e^c = ({\bf 1},
{\bf 1})_1 &
\end{array}
\end{equation}
In addition it includes a Higgs field $h = ({\bf 1}, {\bf
2})_{-1/2}$ which may soon be discovered at the LHC.

At first sight, these groups and matter representations may look
rather random. However it has been known for some time that these
representations fit together very nicely if we embed the gauge group
into a larger semi-simple group \cite{Georgi:1974sy}. The simplest
choice is
\begin{equation}
SU(3) \times SU(2) \times U(1)\ \subset\ SU(5)
\end{equation}
In this case, the matter representations are combined into
\begin{equation}
\bt\ =\ (Q,u^c,e^c), \qquad \bfb\ =\ (L,d^c)
\end{equation}
where $\bt$ denotes the two-index ant-symmetric representation of
$SU(5)$, and $\bfb$ denotes the anti-fundamental representation of
$SU(5)$.

This `unification' of the representations can be continued. The next
step is $SU(5) \subset SO(10)$, which unifies the matter
representations as well as the neutrino into a single spinor
representation of $SO(10)$:
\begin{equation}
{\bf 16}\ =\ \bt + \bfb + {\bf 1}
\end{equation}
where the singlet corresponds to the right-handed neutrino. One may
continue with $SO(10) \subset E_6$, which unifies matter with the
Higgs field.

This remarkable bit of group theory could be a simple mathematical
curiosity, explained to some extent by anomaly
cancellation. However, there is important dynamical evidence that
nature makes use of this idea. One of the most compelling pieces of
evidence is obtained by adding supersymmetric partners around a TeV,
and extrapolating the one-loop running of the gauge couplings
upwards. One arrives at the following remarkable picture
\cite{Georgi:1974yf,Dimopoulos:1981yj}:

\begin{figure}[th]
\begin{center}
            \scalebox{.8}{
  \includegraphics[width=\textwidth]{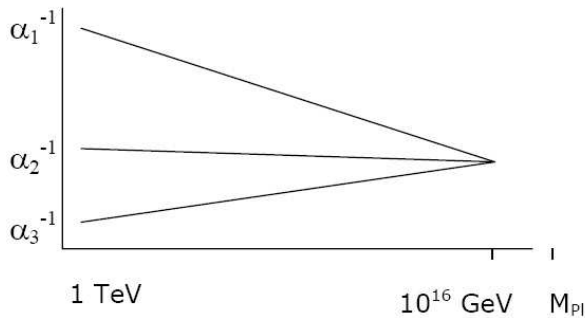}
}
\end{center}
 \caption{Important scales in particle physics. The
weak scale $\sim 1$ TeV is associated with electro-weak symmetry
breaking and will be explored by the LHC. Extrapolating the gauge
couplings of the strong force ($\alpha_3$), weak force ($\alpha_2$),
and hypercharge ($\alpha_1$), we see that they meet at an energy
scale $\sim 10^{16}$ GeV which is called the GUT scale. At the
Planck scale $\sim 10^{19}$ GeV gravity becomes strong and we need a
full theory of quantum gravity.} \label{firstfig}
\end{figure}

The unification of couplings gives strong evidence for the idea that
the three forces of the Standard Model merge into a single force at
high energy scales \cite{Georgi:1974sy}. Such gauge theories go
under the name of Grand Unified Gauge Theories or simply GUTs. The
convergence of couplings isn't just an accident of one number coming
out right. The fact that the couplings meet just below the Planck
scale is also very suggestive. And the picture is fairly robust and
in principle independent of supersymmetry;
for example it is not
materially affected by adding complete $SU(5)$ multiplets to the
theory at intermediate energy scales. Without supersymmetry though,
the couplings don't meet with the same impressive accuracy.

Furthermore, these models not only provide an explanatory framework
for some known patterns, but they also yield a new and very
characteristic prediction. As was emphasized by Georgi and Glashow,
GUT models collect particles with different baryon number in the
same multiplet and therefore lead  to proton decay, so looking for such
decay is one of the best ways to test the hypothesis. In fact, proton decay
is one of the few ways we have to directly probe energy scales as
large as $10^{16}$ GeV. It is one of the outstanding questions of particle
physics and remains on the experimental agenda.
For constraints on Grand Unified models from proton decay, both with and
without supersymmetry, see eg. \cite{Raby:2002wc,Nath:2006ut,Senjanovic:2009kr}.

\subsection{The heterotic string}

$ $

 String theory gives a further generalization of this picture by
also adding quantum gravity into the mix. We would like to
understand how to recover Grand Unification from string theory. This
was first achieved by \cite{Candelas:1985en} in '85.

We start with the ten-dimensional $E_8\times E_8$ heterotic string.
Of course its symmetry group is much too big. We want to break this
to four-dimensional $N=1$ super-Poincar\'e invariance and $SU(5)$
gauge group. To break the Poincar\'e symmetry, we take the space
time of the form
\begin{equation}
{\bf R}^{1,3} \times X_6
\end{equation}
It is possible to make a more general Ansatz with a warp factor, but
we will not consider it as it is still poorly understood. Preserving
$4d$ $N=1$ super-Poincar\'e then implies that $X_6$ must be a
Calabi-Yau three-fold. To break the gauge symmetry we turn on
a non-trivial gauge field along $X_6$. Four-dimensional
supersymmetry then implies that the $E_8\times E_8$ connection must satisfy
the Hermitian Yang-Mills equations:
\begin{equation}\label{HYM}
F^{0,2} = 0, \quad g^{i\bar j} F_{i\bar j} = 0
\end{equation}

Given a solution $(X_6,V)$ to these equations, we may then calculate
the effective four-dimensional gauge theory by computing certain
Dolbeault cohomology groups. The Standard Model gauge group and the
charged matter fields (the `visible sector') all come from
Kaluza-Klein reduction of the $10d$ gauge theory, using just
one of the two $E_8$ groups. They are
given by the Dolbeault cohomology groups
\begin{equation}
H^p(X_6,V_{E_8})
\end{equation}
where $V_{E_8}$ is the adjoint bundle for one of the $E_8$ gauge
groups. The second $E_8$ is said to yield a hidden sector.
For $p=0$ this cohomology group counts generators of the $4d$ gauge group,
and for $p=1$ this counts $4d$ matter fields. So we want to find
pairs $(X_6,V_{E_8})$ such that
\begin{eqnarray}
H^0(X_6,V_{E_8}) &=& su(3) \oplus su(2) \oplus u(1)\\
H^1(X_6,V_{E_8}) &=& \sum_{i=1}^3 (\bt_i \oplus \bfb_i) \oplus higgs
\end{eqnarray}
We used $SU(5)$ notation for the quarks and leptons to emphasize
that they come in complete $SU(5)$ multiplets, even though the GUT group
is broken.
Furthermore, there is a Yoneda product
\begin{equation}
H^1(X_6,V_{E_8})\times H^1(X_6,V_{E_8})\times H^1(X_6,V_{E_8})\ \to\
{\bf C}
\end{equation}
which computes the $4d$ Yukawa couplings. Here we used the Calabi-Yau
condition and the obvious anti-symmetric three-index tensor for
$E_8$ to get a number. The K\"ahler potential
(kinetic terms) may not be computed exactly, however there are
numerical techniques that allow one to approximate it
\cite{Douglas:2006hz}.

Note that even ignoring gravity, string theory has added another
level of unification to the story of Grand Unification. In
conventional $4d$ GUTs, matter and gauge fields still appeared as
separate entities. But by adding extra dimensions we have managed to
unify gauge and matter fields into a single structure: the unique
$10d$ supersymmetric gauge theory with $E_8$ gauge group.

\subsection{The landscape}

$ $

A solution $(X,V)$ gives only a first (tree level) approximation to
the $4d$ physics. Simple solutions tend to have a large number of moduli,
and there are no-go theorems that say that suitable classical
solutions without moduli do not exist. The presence of light moduli
contradicts known experimental facts, such as the equivalence
principle or Big-Bang nucleosynthesis. Furthermore the correct
theory of particle physics below a TeV (i.e. the Standard Model) is
not supersymmetric. Finding a stable weakly coupled solution is a
complicated dynamical problem, and unfortunately many predictions
depend on it.

It is by now appreciated that many of the moduli of $(X,V)$
(in particular complex structure moduli, vector bundle moduli and
even some K\"ahler moduli) can be lifted at tree level, leaving only
very few moduli $S_i$ to be stabilized by quantum effects.
However the number of such classical solutions $(X,V)$ which we would
like to use as a first approximation to a true solution, keeping $S_i$ fixed,
appears to be very large. The
existence of an exponentially large number of
solutions is well-advertized in the context of $F$-theory flux vacua
\cite{Bousso:2000xa,Denef:2004ze}. It is less well-known that there
is an equally large heterotic landscape. This was overlooked for
many years, but should have been expected on
the basis of heterotic/$F$-theory duality. In particular, the construction of
heterotic bundles on elliptic Calabi-Yau manifolds was related in
\cite{Donagi:2009ra} to a Noether-Lefschetz problem, which leads to an
exponentially large number of
solutions, in the same manner as the $F$-theory landscape.

Let us consider the construction of the MSSM in the heterotic string,
on a fixed elliptically fibered Calabi-Yau three-fold $X$.
Apart from the Calabi-Yau, the main ingredient is a stable rank five
bundle $V$  with $c_1(V)=0$ and ${1\over 2}
c_3(V) = 3$. There is a bound on $c_2(V)$ from tadpole cancellation,
which will show up as a cut-off $\Lambda$ below.
Using the results of \cite{Donagi:2009ra} we can make
the following rough, Bousso/Polchinski-like estimate for the number
of isolated solutions:
\begin{equation}
 \#\, {\rm solutions} \ \simeq\ {1\over (r/ 2)!}\,(2\pi  \Lambda)^{\rm r/2} \ \simeq \
\left( {4\pi e \Lambda\over r}\right)^{r/2} \ \simeq \ 10^{1000}
\end{equation}
Here $\Lambda$ is a cut-off due to tadpole cancellation, which we
took conservatively of order $\Lambda \sim r/24$. The number
$r \sim 10^3$ to $10^5$ comes from the rank of the $H_2$ homology lattice of
a degree five spectral cover over ${\bf CP}^2$. This version of the landscape
arises solely from compactifying the
$10d$ $E_8$ gauge theory, i.e. from the choice of solution to the
hermitian-Yang-Mills equation on a rank five bundle, holding
everything else (in particular the background $X$) fixed. As such,
it directly affects the parameters in the visible sector.

These numbers are so large that is simply pointless to try to
enumerate the set of MSSMs, even on a fixed $X$. However
generic solutions are expected to have
qualitatively similar particle physics. This lends support to the
traditional idea of naturalness in model building: unless there is
extra well-motivated structure, dimensionless parameters should be
assumed to be of order unity.

\subsection{Extension to seven, eight and nine dimensions}
\label{789Extension}
$ $

In '85 our picture of string theory was rather limited, and the only
known realization of $E_8$ gauge theory was in the context of the
heterotic string. However during the duality revolution in the '90s
we learned that lower dimensional $E_8$ gauge theories are also
realized in string theory. Thus our excuse for only
constructing GUT models from ten-dimensional gauge theory has
disappeared. These constructions were carried out in a series of
papers in the last three years (some yet to appear).

\begin{figure}[t]
\begin{center}
\renewcommand{\arraystretch}{1.5}
\begin{tabular}{|c|c|}
  \hline
  \ \  {\it dim} \qquad & \qquad  {\it stringy realization}\qquad \qquad \\
  \hline \hline
  10d \qquad & $E_8 \times E_8$ heterotic string \\
  9d & strongly coupled type I' \\
  8d & $F$-theory on ALE \\
  7d & $M$-theory on ALE \\
  \hline
\end{tabular}
\\[5mm]
\parbox{10cm}
{\bf Table 1: \it Branes with exceptional gauge symmetry in string
theory.}
\renewcommand{\arraystretch}{1.0}
\end{center}
\end{figure}

The basic summary for the realization of $E_8$ gauge theories in
string theory is given in table $1$. There are no supersymmetric
gauge theories above ten dimensions, which is why the table ends at
the upper end. It is possible to realize $E_8$ gauge theory below 7
dimensions, however then there is not enough room to get chiral
matter in 4d. This is why the table ends at the lower end.

The entries for
$d=7,8,9$ look slightly mysterious: they do not correspond to weakly coupled
string theories and the way the non-abelian gauge symmetries arise is not
completely obvious, as it involves non-perturbative physics. Obviously we cannot do justice to it here,
but let us at least give a very quick summary.

The starting point is $M$-theory, the conjectural non-perturbative completion of eleven-dimensional
supergravity. $M$-theory
on a smooth space-time does not give rise to non-abelian gauge symmetry, but $M$-theory on a singular ALE space
of type ADE gives rise to seven-dimensional super-Yang-Mills theory localized at the
singularity, with the corresponding ADE gauge group \cite{Witten:1995ex}. An important consistency check
is that when we resolve the singularity, quantized $M2$-membranes wrapping the vanishing cycles
reproduce the off-diagonal components of the Yang-Mills multiplet (i.e. those not proportional to the
Cartan generators). $M$-theory on an elliptically fibered space-time with section can
be mapped to ten-dimensional type IIb supergravity on the section times a circle. In the limit that the elliptic
fiber shrinks to zero, the circle decompactifies.
The modular parameter of the elliptic fiber
is identified with a varying axio-dilaton of the IIb supergravity, $\tau = a + i e^{-\phi}$, and the mechanism
of non-abelian gauge enhancement is similar to $M$-theory. This is called $F$-theory \cite{Vafa:1996xn}.
$M$-theory on an $S^1/Z_2$
orbifold gives rise to $10d$ $E_8$ gauge theory on each of the two walls. This is the Horava-Witten picture
\cite{Horava:1995qa}.
By shrinking the interval, we recover the weakly coupled heterotic string.
By instead compactifying on $S^1$ and shrinking it, we get the type I' theory \cite{Polchinski:1995df}.
In each of these contexts, when the non-abelian
gauge symmetry is of type
$SU(n)$, $SO(n)$ or $USp(n)$, we can often extrapolate to a weakly coupled $D$-brane configuration in
perturbative string theory, but for the exceptional
cases this is not possible. These issues were only understood after the non-perturbative behaviour of string theory became clearer
in the '90s.

In this review we will be
firmly focused on the Yang-Mills theory itself, without much consideration of its origins.
Part of the reason for this is that if we want to consider exceptional gauge groups in higher dimensions, then
the string point of view involves singular geometries and strongly coupled physics, whereas the higher dimensional
Yang-Mills theory is weakly coupled in the derivative expansion.
We will briefly explain how to connect with string theory in section \ref{SpectralALE}.
The upshot will be that we can reconstruct the local geometry from the Yang-Mills theory.
The type I' story requires a separate discussion, which will not be attempted here.

Unification with eight-dimensional Yang-Mills theory was initiated in
\cite{Donagi:2008ca,Beasley:2008dc,Hayashi:2008ba} and is currently
the most developed and best understood of the new cases. The main reason for this is
that in even dimensions one can use the powerful techniques of
complex and algebraic geometry. The seven-dimensional story was
initiated in \cite{Pantev:2009de} and the nine-dimensional story
will appear in \cite{PW2}.

The strategy that we will use is similar to the heterotic string.
The main new idea is that instead of Hermitian-Yang-Mills bundles,
we have to consider their close cousins: Higgs bundles. This
consists of a bundle $E$ together with a map
\begin{equation}
 \phi: E \to E \otimes N
\end{equation}
where $N$ is a suitable twisting bundle. The connection $A_\mu$ on
$E$ and the `Higgs field' $\phi$ have to satisfy certain versions of
Hitchin's equations, which we discuss below. (The terminology is
standard and the field $\phi$ has no direct relation to the Higgs
boson $h$ of the Standard Model).

\section{Higgs bundles}

\subsection{Dimensional reduction}

$ $

The maximally supersymmetric Yang-Mills theory in $p$ dimensions is
uniquely determined by the choice of gauge group. Its Lagrangian can
be obtained by dimensional reduction from the $10d$ supersymmetric
Yang-Mills theory. Therefore we would expect that the BPS conditions
can be obtained from dimensional reduction of the $6d$ hermitian
Yang-Mills equations (\ref{HYM}).

There is a well-known analogous story involving reduction of the
$4d$ hermitian Yang-Mills equations, aka the ASD equations.
Reducing to three dimensions yields the Bogomolnyi equations:
\begin{equation}
F = * D_A \phi
\end{equation}
The reduction to two dimensions is known as Hitchin's equations:
\begin{equation}\label{HitchinEquation}
\bar \partial_A \phi = 0, \quad \Lambda F + [\phi,\phi^*] = 0
\end{equation}
where $\Lambda$ denotes contraction with the K\"ahler form. These
equations turn out to have many applications, and several variants
have been considered. We will be interested in some of its
higher dimensional variants.

Indeed let us now consider the dimensional reduction of the $6d$
hermitian Yang-Mills equations (\ref{HYM}), relevant for
compactifications of $10d$ Yang-Mills down to four dimensions. For
$8d$ Yang-Mills on ${\bf R}^{1,3} \times S_4$ we get $F^{0,2} = 0$,
as well as equations (\ref{HitchinEquation}) above on $S_4$. For
$7d$ Yang-Mills on ${\bf R}^{1,3} \times Q_3$ we get
\begin{equation}\label{MHiggs}
F - [\phi,\phi] = 0, \quad D_A\phi = 0, \quad D_A^\dagger \phi = 0
\end{equation}
This is the real version of Hitchin's equations, defined on a real
three-manifold $Q_3$. The gauge and Higgs field can be assembled into a complex gauge field
$A+i\phi$, and the first two equations above assert that this is a
complex flat connection on $Q_3$.
For $9d$ Yang-Mills we compactify on ${\bf
R}^{1,3} \times X_5$, where $X_5$ is a five-manifold given by an
$S^1$-fibration over a K\"ahler surface. Denote complex coordinates
on the base by $z$ and along the $S^1$ by $y$. Then we get
\begin{equation}\label{IpHiggs}
F_\perp^{0,2}=0, \quad F_{\bar z y} = i D_{\bar z}\phi, \quad
\Lambda F_\perp + i D_y \phi = 0
\end{equation}
where $\perp$ denotes the components orthogonal to the $S^1$. This
is a higher dimensional generalization of the Bogomolnyi equations.

\subsection{Compactification of eight-dimensional SYM}

$ $

We can now summarize the main results of the Kaluza-Klein reduction
of the higher dimensional Yang-Mills theory to four dimensions. We
take the $8d$ supersymmetric Yang-Mills theory as our main example
\cite{Donagi:2008ca,Beasley:2008dc,Donagi:2009ra,Hayashi:2009ge}. The
Higgs bundle will be defined on a compact K\"ahler surface $S$. As
we have seen, the bosonic fields are given by a connection $A_\mu$
on a bundle $ad(E)$ on $S$, and a complex Higgs field $\Phi$ valued
in $ad(E) \otimes K_S$. The twisting by the canonical bundle is
required to preserve the $4d$ super-Poincar\'e invariance. The equations
are given by
\begin{equation}
F^{0,2} = 0, \quad \bar{\partial}_A\Phi = 0, \quad \Lambda F +
[\Phi,\Phi^*] = 0
\end{equation}
In practice, we actually need meromorphic Higgs bundles, or in recent physics language,
we need to insert certain surface operators along a curve $\Sigma_\infty$ in $S$.
The simplest such defects lead to parabolic Higgs bundles \cite{Donagi:2011jy}.
They may be thought of as arising from
integrating out hypermultiplets living on $\Sigma_\infty$, which become dynamical when
we embed the gauge theory into a compact model. We will
not discuss these defects explicitly here.

Given a solution to these equations, we would like to know the
effective gauge theory in four dimensions obtained from the
Kaluza-Klein reduction of the eight-dimensional gauge theory on
${\bf R}^{1,3} \times S$. We define a two-term complex ${\mathcal
E}^\bullet$:
\begin{equation}
{\mathcal E}^\bullet\ = \ ad(E)\ \mathop{\to}^{\Phi}\ ad(E) \otimes
K_S
\end{equation}
The $4d$ gauge theory is derived by computing the hypercohomology
groups
\begin{equation}
{\mathbb H}^p(S,{\mathcal E}^\bullet)
\end{equation}
For $p=0$ this counts generators of the $4d$ gauge group, and for
$p=1$ it counts $4d$ matter fields. We can compute the Yukawa
couplings from a Yoneda product
\begin{equation} {\mathbb H}^1(S,{\mathcal E}^\bullet) \times {\mathbb
H}^1(S,{\mathcal E}^\bullet) \ \to \ {\mathbb H}^2(S,{\mathcal
E}^\bullet)
\end{equation}
and conjecturally we can even numerically approximate the
hermitian-Einstein metric and therefore the kinetic terms
\cite{Donagi:2011jy}. Of course this parallels the analogous
statements for the heterotic string we discussed earlier.

\begin{figure}[t]
\begin{center}
            \scalebox{.6}{
               \includegraphics[width=\textwidth]{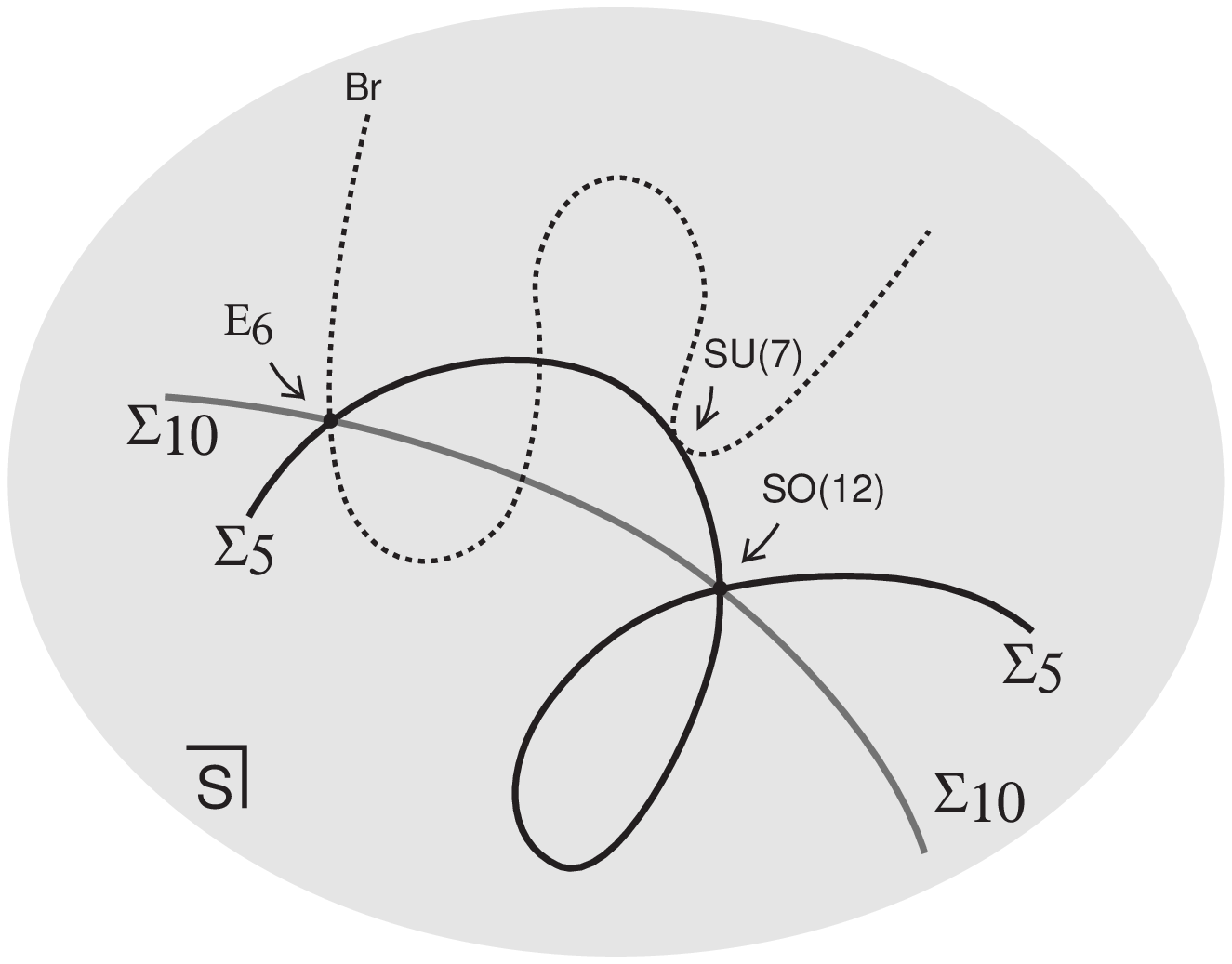}
               }
\end{center}
\caption{ \it Schematic picture of the matter curves, branch locus
and their intersections in an $SU(5)$ GUT model, for generic values
of the complex structure moduli.}
\label{Enhancements}
\end{figure}

These statements perhaps seem somewhat abstract. It is often
possible to give more intuitive pictures for the wave functions of
the massless modes. The main picture in this regard is that a
generic zero of $\Phi$ can be interpreted as a vortex string on $S$, and as
is well known one tends to get charged zero modes localized on such
a vortex. Thus while the wave functions of the four-dimensional
gauge fields are spread out over $S$, the wave functions of charged
matter such as the $\bt$ or $\bfb$ of an $SU(5)$ GUT tends to be
localized on some Riemann surface in $S$, see figure
\ref{Enhancements}. For $7d$ super-Yang-Mills theory compactified on $Q_3$
we get a similar picture, with gauge fields spread over $Q_3$ but
chiral matter in the $\bt$ or $\bfb$ generically localized at points on $Q_3$.
Such localization leads to interesting
possibilities for phenomenology. But for non-generic configurations
the intuition can be misleading,
and the advantage of the above formulation is that it is precise and
general.

%

\subsection{Higgs/spectral correspondence}
\label{HiggsSpectral}
$ $

Now that we have the general statements, we need a good method for
constructing solutions and computing the hypercohomology groups.
Hitchin's equations are highly non-linear, so trying to solve them
in closed form is bound to fail. Instead we will make use of two
standard techniques. The first is the technique of splitting up the
equations in `$F$-terms' and `$D$-terms', i.e.:
\begin{itemize}
  \item a pair of complex equations: $F^{0,2} = \bar \partial_A \Phi = 0$\\
  \item a moment map: $\Lambda F + [\Phi,\Phi^*] = 0$
\end{itemize}
The strategy is then to temporarily ignore the $D$-terms, and focus
on the $F$-terms, which we can solve exactly up to {\it complexified}
gauge transformations. This amounts to focussing on the holomorphic structure
and ignoring the hermitian metric. Once we have done
that, we can go back and try to apply some powerful results in
geometric invariant theory regarding $D$-flatness.

To solve the $F$-term equations, we use the second essential tool,
namely the Higgs/spectral correspondence. Essentially it turns the problem
of constructing solutions to the $F$-terms into an `abelianized'
version which is much simpler to solve. To avoid the somewhat
intricate group theory of $E_8$, we will take $E$ to be a rank $n$
vector bundle transforming in the fundamental representation of
$U(n)$.

The Higgs/spectral correspondence states that there is an
equivalence
\begin{equation}
{\rm Higgs\ bundle\ }(E,\Phi)\  \leftrightarrow\ {\rm spectral\
sheaf\ } \mathcal{L}
\end{equation}
mapping the (holomorphic data of the) Higgs bundle to its spectral
data. That is, we interpret $\Phi$ as a map $E \to E \otimes K_S$,
and then we fiberwise replace $\Phi$ by the eigenvalues and $E$ by
the eigenvectors. The relation can concisely expressed through a
short exact sequence
\begin{equation}
0\ \to\ \pi^*E\ \mathop{\longrightarrow}^{\Psi}\ \pi^*(E\otimes
K_S)\ \to\ {\mathcal L}\ \to\ 0
\end{equation}
where $\Psi \equiv \Phi - \lambda I_{n\times n}$, $\pi^*$ is
pull-back to the total space of the canonical bundle, and $\lambda$
is the tautological section of $\pi^*K_S$. For future use, we also
denote the total space of $K_S$ by $X$. The similarity to our notation
for the heterotic Calabi-Yau three-fold is not entirely coincidental.
There are versions of this correspondence also for meromorphic Higgs bundles.

Such spectral sheaves are fairly simple to write down explicitly,
certainly compared to solving non-abelian equations. As
an example relevant for phenomenological models, let us consider breaking an $E_8$ gauge group to
$SU(5)_G$. This requires writing down an $Sl(5,{\bf C})_H$ Higgs
bundle over $S$, where the structure group $H$ is the commutant of $G$ in the
(complexified) $E_8$.
By the Higgs/spectral correspondence, the data of such an $Sl(5,{\bf C})_H$ Higgs
bundle is given by
\begin{itemize}
  \item the spectral cover $C$ in $X$ (the support of $\mathcal{L}$), given explicitly by a
  degree five polynomial;\\
  \item the spectral line bundle $L$ in ${\rm Pic}(C)$, which may be specified by
  writing an explicit divisor $D$ on $C$ and putting $L = {\mathcal O}(D)$. We then define
  ${\mathcal L} = i_*L$, where $i$ is the inclusion $C \hookrightarrow X$ and $i_*$ is the push-forward.
\end{itemize}
Here we can see the Noether-Lefschetz problem rearing its head. We
need to choose a curve in a surface in a three-fold. Simple
examples constructed by hand will have many moduli. But if the curve is
sufficiently `generic' then deforming the surface or the three-fold
will destroy the curve, in other words the corresponding deformation moduli are
stabilized. However finding the most rigid solutions requires us to solve a complicated
Noether-Lefschetz problem, and the number of solutions grows exponentially
with the rank of $H^2(C)$. This brings us back to the landscape problem. If these vacua
are really stable after  including quantum corrections (and most are probably not), then
we'd have to ask who lives in the
other mathematically consistent universes, and how was our universe able to solve the
complicated computational problem of finding the most rigid, stable solutions. It is probably
best to focus on general features of the whole class of solutions, rather than on
individual solutions.

Leaving such questions aside, another important point about the separation into
$F$-terms and $D$-terms is that changing the hermitian metric on $E$
does not affect the hypercohomology groups ${\mathbb
H}^p(S,{\mathcal E}^\bullet)$. Therefore the low energy spectrum and
Yukawa couplings are independent of the hermitian metric on $E$, and
in fact they may be computed directly in terms of the spectral sheaf
$\mathcal{L}$ using Ext groups. Needless to say, this simplifies
one's life tremendously.

This takes care of the $F$-terms. The $D$-terms (i.e. the moment map equation)
may now be
interpreted as an equation for the hermitian metric on $E$, and a
solution (if it exists) is called a hermitian-Einstein metric. As
mentioned already, solving this equation in closed form is virtually
impossible. Instead one would like to appeal to a Higgs bundle
version of the Donaldson-Uhlenbeck-Yau theorem:
\begin{center}
\begin{quotation}
A unique HE metric exists $\ \leftrightarrow\ $ the Higgs bundle is
poly-stable.
\end{quotation}
\end{center}
A Higgs bundle is stable if for every Higgs sub-bundle $F \subset E$, we have
\begin{equation}
\mu(F) < \mu(E)
\end{equation}
The slope $\mu(E)$ is defined to be the ratio of the degree of $E$
(with respect to some choice of K\"ahler form) and the rank of $E$.
A Higgs bundle is said to be poly-stable if it is the sum of stable
Higgs bundles with the same slope.

The beauty of this type of statement is that existence of the
hermitian-Einstein metric is a differential-geometric criterion,
whereas poly-stability is an algebro-geometric
criterion. In particular, one can easily translate the requirement of stability for the Higgs bundle
$(E,\Phi)$ to stability for the spectral sheaf
$\mathcal{L}$. Unfortunately the proven versions of this
correspondence do not suffice in our context, because we really need
certain meromorphic incarnations of Higgs bundles. In particular, the
stability condition depends on boundary data associated to the
defect/surface operator, as is well-known for parabolic Higgs bundles.
Modulo that issue however, we see that the Higgs/spectral correspondence
answers the problem posed at the beginning of this section:
constructing solutions and understanding the $4d$ effective theory boils
down almost entirely to constructing and studying suitable spectral
sheaves $\mathcal{L}$.

In other dimensions, the story is in principle similar, but less understood.
In $M$-theory we were dealing with (\ref{MHiggs}), which describes complex
flat connections on a real manifold $Q_3$. Up to complexified gauge transformations,
we can exchange such a complex flat
connection for spectral data. In this case, the spectral cover corresponds to a Lagrangian $A$-brane
in the cotangent bundle $T^*Q_3$ \cite{Pantev:2009de}, and the $4d$ spectrum is computed
by Floer cohomology groups of the $A$-brane.

\subsection{Spectral/ALE correspondence}
\label{SpectralALE}
$ $

Although we have claimed that our constructions are realized in
string theory, the reader may be puzzled that the descriptions we
have obtained seem to look very different from the traditional
descriptions of $F$-theory, $M$-theory and type I', which we briefly
touched on in section \ref{789Extension}. Four-dimensional
compactifications of $F$-theory
are usually said to involve elliptically fibered
Calabi-Yau manifolds with a configuration for a certain three-form
(mathematically, a 2-gerbe with some additional properties).
Four-dimensional compactifications of $M$-theory are usually said
to involve manifolds with $G_2$-holonomy and a flat three-form
(mathematically, a flat 2-gerbe). So far these did not yet
appear in our story. In order to see the relation, we need another
dictionary: the spectral/ALE correspondence
\cite{Curio:1998bva,Donagi:2009ra,Pantev:2009de}.

In order to avoid the somewhat complicated representation theory of
$E_8$, we will illustrate the correspondence by focussing on $U(n)$
Higgs bundles and $A_{n-1}$-type ALE fibrations. Consider the ALE
fibration $Y \to S$ given by
\begin{equation}
y^2\ =\ x^2 + a_0 z^n + \ldots + a_n
\end{equation}
The $a_i$ are sections of certain line bundles on $S$. Let us define
a fibration $R \to S$ whose fibers over a point on $S$ are the $n$
lines defined by
\begin{equation}
y\ =\ x, \quad  a_0 z^n + \ldots + a_n\ =\ 0
\end{equation}
The second equation specifies $n$ points in the $z$-plane, hence
defines an $n$-fold covering $C \to S$, which we identify with the
spectral cover. We also get a natural projection $p : R \to C$ by
replacing each line with a point. We further have a natural
inclusion $i: R \hookrightarrow Y$. Then the spectral/ALE
correspondence is given by the `cylinder' maps $i_*p^*$ and
$p_*i^*$:
\begin{equation}
H^{i,j}(\bar C) \  \leftrightarrow \  H_p^{i+1,j+1}(\bar Y)
\end{equation}
where $H_p$ denotes a certain primitive part of the cohomology, and
$\bar C$ and $\bar Y$ denote certain compactifications. In
particular this maps the class of the spectral line bundle in ${\rm
Pic}(\bar C)$ to the Deligne cohomology class of a $2$-gerbe on
$\bar Y$. The more sophisticated version for exceptional gauge
groups was described in \cite{Curio:1998bva,Donagi:2009ra}.
Essentially this mapping establishes a isomorphism of certain Hodge
structures associated to the spectral cover and ALE sides.

The spectral/ALE correspondence exists for any gauge group, but the group theory
tends to get more involved when we go away from the unitary groups.
It so happens that for $SU(5)_{GUT}$ model building,
which is the most relevant case phenomenologically, the group theoretic aspects
simplify somewhat. A configuration in an $E_8$ gauge theory with an unbroken
$SU(5)_{GUT}$ group corresponds to the following ALE fibration over $S$ \cite{Donagi:2009ra}:
\begin{equation}
\label{E8Unfolding}
y^2 \ = \ x^3 +  a_0 z^5 + a_2 z^3 x + a_3 z^2 y + a_4 z x^2 + a_5 xy
\end{equation}
Here the $a_i$ are sections of certain line bundles over $S$.
Note that we can easily put this in the form of a Weierstrass
model by making some coordinate redefinitions. Under the
spectral/ALE correspondence this gets mapped to an $E_8$
spectral cover which decomposes into various pieces.
One of these pieces is a five-fold cover $C_5$ of $S$, given by a degree
five equation in the total space
of $K_S \to S$:
\begin{equation}
0 \ = \ a_0 s^5 + a_2 s^3 + a_3 s^2 + a_4 s + a_5
\end{equation}
Here $s$ is a local coordinate on the fiber of $K_S$. This is
precisely the spectral cover for the
fundamental representation of the $Sl(5,{\bf C})_H$ Higgs
bundle mentioned in section \ref{HiggsSpectral}.
The intersection of $C_5$ with the zero section
yields a curve $a_5=0$ where the $A_4$ singularity of the
generic ALE fiber (\ref{E8Unfolding}) further degenerates to
type $D_5$, as one can easily check. This is the curve $\Sigma_\bt$
where matter fields on the
$\bt$ and $\btb$ of $SU(5)_{GUT}$ propagate, see figure
\ref{Enhancements}. One also gets a ten-fold covering
$C_{10}$, the spectral cover in the anti-symmetric representation of our
$Sl(5,{\bf C})_H$ Higgs
bundle in section \ref{HiggsSpectral}. Its intersection with
the zero section yields a curve
$a_0 a_5^2 - a_2 a_3 a_5 + a_3^2 a_4 = 0$ where the $A_4$
singularity of the ALE fiber degenerates to $A_5$.
This is the curve $\Sigma_\bfv$ where the $\bfv$ and $\bfb$
matter of $SU(5)_{GUT}$ propagate, again see figure \ref{Enhancements}.
One gets various further degenerations
in higher codimension. For details of this construction, see \cite{Donagi:2009ra}.

Using the spectral/ALE correspondence for $E_8$ we obtain an elliptic Calabi-Yau
four-fold $Y$ with boundary together with a Deligne cohomology
class. We can then proceed to glue $Y$ into a compact Calabi-Yau
four-fold \cite{Donagi:2009ra}. In this way we recover the
traditional description of $F$-theory vacua.

In the $M$-theory context, the spectral data
was given by a Lagrangian $A$-brane
in $T^*Q_3$. Under the spectral/ALE correspondence, it gets mapped to an ALE fibration over $Q_3$
with a flat 2-gerbe,  and with singularities at the location of non-abelian
gauge groups and charged chiral matter. We have argued that this seven-dimensional
non-compact manifold should admit a metric with
$G_2$ holonomy \cite{Pantev:2009de}. This should then be glued into a compact $G_2$ manifold.

Note that in the traditional descriptions, we are dealing
with singular spaces. For example to get an $SU(5)$ gauge theory in
four dimensions from $F$-theory, we needed $A_4$ type quotient singularities along a
section isomorphic with $S$ which further degenerate in higher
codimension. Of course, the problem with singularities is that they are singular.
In order to do physics, we need some way to `smooth'
these singularities. The traditional way to do this is by making a
crepant resolution (in physics terms, using $M$/$F$-duality and
moving out on the Coulomb branch). Then one can quantize solitons
wrapped on the exceptional cycles to deduce some of the physics.
Unfortunately this extrapolation is not valid at the level of
$D$-terms and it obscures many aspects of the physics which are
relevant for phenomenology. For example a proper definition of the
$F$-theory 2-gerbe and the analogue of its hypercohomology groups
has not yet been given, in part because the 2-gerbe may obstruct
such resolutions. Furthermore it is not clear how to extend this
approach to other setting like $M$-theory on $G_2$ manifolds with
singularities, where no natural resolution is available. The Higgs
bundle/Yang-Mills theory approach yields another way to smooth the singularities, which has
proven more useful for the questions discussed here.

\subsection{Some open problems}

$ $

Here we collect some open mathematical problems. They were mostly
already pointed out in the text, but we group them here for
convenience. It is not meant to be a complete list, and we only list
questions that fit in the scope of this review.

\begin{itemize}
  \item {\bf Explicit constructions.} It is an open problem to find non-abelian solutions to the
  $F$-term part of the Yang-Mills-Higgs equations in odd dimensions, given in equations
  (\ref{MHiggs}) and (\ref{IpHiggs}). Furthermore, given a solution, one would like an effective
  method for computing hypercohomology groups.
  For solving (\ref{MHiggs}) one could appeal to Donaldson-Corlette.
  However as we mentioned, for phenomenological applications one should add
  source terms to (\ref{MHiggs}) and (\ref{IpHiggs}) along a submanifold. An interesting sub-problem
  is to classify the most general supersymmetric boundary conditions (source terms).

  A class of non-abelian solutions for the type I' case
  will appear in \cite{PW2}. The strategy there is to assume an $S^1$ symmetry and
  make use of the Fourier-Mukai transform to turn the question into a
  holomorphic problem. For the $M$-theory models we could use a similar strategy, by using mirror
  symmetry for $T^*Q_3$. Presumably there are other approaches.
  Is it possible to give a more general construction?

  The equations encountered in this review often reappear in other contexts, see eg.
  \cite{Kapustin:2006pk} or the literature on $Sl(N)$ knot invariants. The source terms
  are said to be associated to `defect operators.' Thus a positive answer could also be
  helpful in other contexts. \\
  \item {\bf D-flatness.} Once we can construct solutions to the $F$-term equations,
  we still have to solve the moment map equation. This is called proving $D$-flatness
  in physics terminology. The standard strategy is to formulate some type of stability
  condition, and then to try and prove an analogue of the Donaldson-Uhlenbeck-Yau theorem
  for Hermitian-Yang-Mills connections.

  In the holomorphic case one can make natural conjectures for the stability criterion
  (including source terms), but a proof
  of $D$-flatness is missing. Mochizuki has given proofs for certain classes of parabolic Higgs bundles,
  which unfortunately are slightly different from the Higgs bundles considered here. For the odd
  dimensional cases, it is not clear to us what the correct
  stability condition should be.\\
  \item {\bf Deformation theory.} To find the
  low energy effective action, we needed to
  understand the infinitesimal deformations of the Higgs bundle.
  They were classified by certain hypercohomology groups of the Higgs bundle.
  In the spectral cover picture,
  one is asking for the infinitesimal deformations of a coherent sheaf or a
  Lagrangian brane with flat connection. These are classified
  by Ext groups or Floer cohomology groups respectively.

  But we also had a third picture: ALE fibrations with a 2-gerbe, satisfying
  some additional conditions. (The conditions in $F$-theory are somewhat
  complicated to state; see appendix C of \cite{Donagi:2008ca} where a
  compactified version of ALE fibrations is considered). So we have a natural question:
  what classifies the first order infinitesimal deformations in the
  ALE fibration picture? And what are the analogues of the Yoneda pairing and the
  higher Massey-like products, which compute Yukawa couplings and higher order interactions?

  It is not clear to us how to answer this question, and maybe a good answer is not possible. There are many
  well-understood spectral cover configurations that get mapped to
  poorly understood ALE-fibrations. As a simple example, one may consider a non-abelian
  bundle on a degenerate non-reduced cover. This should get mapped to an ALE fibration
  with singularities, and some kind of non-abelian 2-gerbe along the singularities, whatever that means exactly.
  Such configurations cannot be lifted to the resolution, but should be included in the
  general formulation. More complicated examples
  can be found in \cite{Donagi:2011jy}. Thus understanding the
  $2$-gerbe is an important part of the problem. One could further ask about stability conditions
  in this picture.
  This is also poorly understood.\\
  \item {\bf Local versus global.} ALE fibrations are said to be local because
  they are non-compact. For string phenomenology, we think of them as a local piece of a
  compact manifold, and
  eventually we would like to embed them in a concrete global model, i.e. in a compact manifold. In the
  context of $F$-theory there has been a lot of progress on this question. Questions remain, but this is outside
  the scope of this review.
  In the context of $M$-theory, we would like to embed our ALE fibration in a compact
  $G_2$ holonomy manifold. The problem here is obvious: there are very few constructions of compact $G_2$ manifolds
  and none appears to be suitable for compactifying our ALE fibrations. \\
  \item {\bf Finite or infinite landscape.} We discussed some crude
  estimates for the number of solutions of the hermitian-Yang-Mills
  or Hitchin equations which reproduce the supersymmetric standard
  model or a unified extension thereof.
  In toy models of flux vacua it appears that the true
  number of solutions can be much larger. In fact, there is a basic
  question if the true number of solutions is even finite. The
  prevailing
  opinion is that it should be finite, but we
  recently constructed an infinite sequence which
  appears to be stable and evades the known
  no-go theorems  \cite{DW7}. (It should be remembered that proving stability
  is of course the hardest part). Such infinite sequences
  would address some concerns raised in \cite{Dine:1985he}.

\end{itemize}

There are of course also many interesting physical questions. In
particular one would like to get an interesting physical prediction
that distinguishes the extra-dimensional models from conventional
four-dimensional models. Once we have done the Kaluza-Klein
reduction, we are working within four-dimensional effective field
theory. So are there processes where we might see the extra
dimensions in the foreseeable future? The answer appears to be yes,
at least in principle. Proton decay probes the GUT scale and could
provide such characteristic signatures
\cite{Friedmann:2002ty,Klebanov:2003my,Donagi:2008kj}. The extra
dimensions become important at the GUT scale and the proton feels
the higher dimensional interactions. For further questions and developments, see the talk by
S. Sch\"afer-Nameki at this conference \cite{SSN}, and the additional reviews \cite{Heckman:2008rb,Wijnholt:2008db,Heckman:2010bq,Wijnholt:2010zz,Weigand:2010wm}.

{\it Acknowledgements}: I would like to thank the organisers of String-Math 2011 for an 
inspiring conference. I would also like to thank R.~Donagi and
T.~Pantev for collaboration on the material presented here, and R.~Donagi
and the referee
for comments on the manuscript.

\bibliographystyle{amsalpha}

\end{document}